\documentclass[conference]{IEEEtran}
\IEEEoverridecommandlockouts
\usepackage{cite}
\usepackage{amsmath,amssymb,amsfonts}
\usepackage{algorithmic}
\usepackage{graphicx}
\usepackage{textcomp}
\usepackage{xcolor}
\def\BibTeX{{\rm B\kern-.05em{\sc i\kern-.025em b}\kern-.08em
    T\kern-.1667em\lower.7ex\hbox{E}\kern-.125emX}}
\begin{document}

\title{Prime Counting Function identity\\
}

\author{
\IEEEauthorblockN{Suyash Garg}
\IEEEauthorblockA{\textit{Department of Mathematics} \\
\textit{BITS-Pilani Hyderabad Campus}\\
Hyderabad, India \\
f20180594@hydereabad.bits-pilani.ac.in}

}

\maketitle

\begin{abstract}
In this paper, a new formula for $\pi^{(2)} (N)$ is formulated, it is a function that counts the number of semi-primes not exceeding a given number N. A semi-prime is a natural number that is the product of exactly two prime numbers, the two primes in the product may equal each other. Semi-prime numbers are also a case of almost primes. Since a formula for this is already known, a new identity that uses the prime counting function is created by equating the two functions.
\end{abstract}


\section{Introduction}
A semi-prime, also called a 2-almost prime, biprime\cite{b1}\cite{b2}, or $pq$-number, is a composite number that is the product of two (possibly equal) primes. The first few are 4,6,9,10,14,15,21,22...(OEIS A001358). The first few semi-primes whose factors are distinct are 6, 10, 14, 15, 21, 22, 26...(OEIS A006881).
A formula for the number of semi-primes less than or equal to $N$ is given by \cite{b3}
\begin{equation}
    \pi^{(2)}(N)=\sum_{k=1}^{\pi(\sqrt{N})}[\pi(N/p_k)-k+1],
\end{equation}
here, $p_k$ is the $k^{th}$ prime, and $\pi(N)$ is the prime counting function, which is a function that gives the number of prime numbers not exceeding the given number N\cite{b4}.
\\ In this paper, a new way to calculate the number of semi-primes not exceeding $N$, is explored.\\
\section{Semi-Prime Counting function}

\textbf{Lemma 1.}
The number of numbers of the form $2p$ not exceeding a given number $N$, is $\pi (N/2)$, where $p$ is any prime number and $\pi (x)$ is the prime counting function.

\textbf{Proof.} The number of primes such that $2p \le N$ is to be calculated. Just by rearranging, it can be shown that $p \le N/2$. This means p can be any prime smaller than $N/2$, which equals $\pi(N/2)$.

\textbf{Corollary 1.}
The number of numbers of the form $pq$ not exceeding a given number $N$, is $\pi (N/q)$, where $p$ is a given constant prime number, $q$ can be any prime number, and $\pi (x)$ is the prime counting.

\textbf{Proof.}
If $p$ is equal to $2$ then, It would be the same as Lemma 1, So, for proving the corollary, it would be sufficient to replace $2$ with a given prime $p$ in Lemma 1.

\textbf{Lemma 2.}
The number of pairs of the form $(p,q)$ where $p$ and $q$ can be any prime such that $pq \le N$ for a given $N$ is 
\begin{equation}
\sum_{k=1}^{\pi (N/2)}[\pi (N/p_k)]    
\end{equation}

\textbf{Proof.}
In Corollary 1, each prime $p_k$ less than N can be used as the constant prime $p$,\\ for $p=2$, number of pairs of the form $(2,q)= \pi (N/2)$,\\ for $p=3$, number of pairs of the form $(3,q)= \pi (N/3)$,\\for $p=p_k$, number of pairs of the form $(p_k,q)= \pi (N/p_k)$,\\ Now, if all of these were to be added, the summation acquired would be, $$\sum_{k=1}^{\pi (N)}[\pi (N/p_k)]$$ Since, all the primes till $N$ were considered for $p$, and $q$ can be any prime, therefore the above formula suffices as a formula for the number of pairs of the form $(p,q)$ such that $ pq \le N$.\\ Also Since, for all primes $p>N/2$, $N/p < 2$ , which means $\pi (N/p) = 0$, therefore the upper limit of the summation can be reduced to $\pi (N/2)$. So, the new formula will be $$\sum_{k=1}^{\pi (N/2)}[\pi (N/p_k)]$$

\textbf{Note.} In Lemma 2, for all the pairs $(p,q)$, where $p$ and $q$ are distinct primes, if order doesn't matter, then each pair is counted twice.

\textbf{Lemma 3.}
The number of primes such that $p^2 \le N$, where $p$ can be any prime is $\pi (\sqrt{N})$.

\textbf{Proof.}
Since $p>0$ and $N>0$,we can get the inequality $p \le \sqrt{N}$.The number of primes such that $p^2 \le N$ is equal to the number of primes such that $p \le \sqrt{N}$, which is equal to $\pi(\sqrt{N})$.

\textbf{Theorem 1.}
\begin{equation}
    \pi^{(2)}(N)=((\sum_{k=1}^{\pi((N/2))}[\pi(N/p_k)]) + \pi (\sqrt{N}))/2,
\end{equation}
is the function that gives the number of semi-primes not exceeding a given number $N$, i.e. a semi-prime counting function.

\textbf{Proof.}
Lemma 2 takes into account all the number of pairs $(p,q)$, such that $p$ and $q$ can be any primes, and $(p,q)$ and $(q,p)$ are considered as $2$ distinct pairs. So, It is trivial that, all the pairs where $p$ and $q$ are distinct will be counted twice and pairs of the form $(p,p)$ will be counted  once. The number of pairs of the form $(p,p)$, is equivalent to the number of primes such that $p^2 \le N$, and Lemma 3 helps get the answer to this question. Now, if this is added to (2), the new formula would be
\begin{equation}
    \sum_{k=1}^{\pi (N/2)}[\pi (N/p_k)] + \pi (\sqrt{N})
\end{equation} In (4) each and every pair is counted twice. so if (4) is halved, the number of pairs of the form $(p,q)$ such that $p$ and $q$ can be any primes, without considering order, can be obtained. This formula will be equivalent to the formula for semi-primes less than or equal to $N$, as, if we multiply the numbers $p$ and $q$ in all the given pairs, we get a list of all the semi-prime less than or equal to $N$, and since we have counted each such pair only once, there would be no repetitions in this list. So, counting the number of elements in this list is equivalent to counting the number of pairs as described above. Therefore, the final formula for the semi-prime counting function is, $$((\sum_{k=1}^{\pi((N/2))}[\pi(N/p_k)]) + \pi (\sqrt{N}))/2$$\\
\section{Example}
Now, there are two formulas for semi-prime counting function.\\

\textbf{Example 1.}
\\Let's take $n=25$.\\Semi-primes not exceeding $25$ are\\$4,6,9,10,14,15,21,22,25$\\
so, 
\begin{equation}
    \pi ^{(2)}(25)=9.
\end{equation}
Using (1)\\
$$\pi ^{(2)}(25)=\sum_{k=1}^{\pi(5)}[\pi(25/p_k)-k+1]$$
$$\pi (5)=3, \pi (8.33)=4, \pi (12.5)=5$$
$$\pi ^{(2)}(25)=(5-1+1)+(4-2+1)+(3-3+1)$$
\begin{equation}
\pi ^{(2)}(25)=5+3+1=9    
\end{equation}
Using (3)\\
$$\pi^{(2)}(25)=((\sum_{k=1}^{\pi((12.5))}[\pi(25/p_k)]) + \pi (5))/2$$
$$\pi (12.5)=5, \pi (8.33)=4, \pi (5)=3, \pi(3.57)=2, \pi(2.27)=1$$
$$\pi^{(2)}(25)=((5+4+3+2+1)+3)/2$$
\begin{equation}
    \pi^{(2)}(25)=(15+3)/2=18/2=9
\end{equation}
$$ (5)=(6)=(7)$$
Therefore, both formulas are correct.
\section*{Identity}
\textbf{Theorem 2.}
\begin{equation}
    \sum_{k=1}^{\pi(\sqrt{N})}[\pi(N/p_k)]-\sum_{k=\pi(\sqrt{N})+1}^{\pi(N/2)}[\pi(N/p_k)]=(\pi(\sqrt{n}))^2
\end{equation}
Where, $\pi (x)$ is the prime counting function, $p_k$ is the $k^{th}$ prime.

\textbf{Proof.}
As proved in theorem 1 and illustrated with Example 1, both of the formulas give the correct count of the number of semi-primes not exceeding a given number $N$, it is logical to equate the two formulas, i.e. $$(1)=(3)$$ $$\sum_{k=1}^{\pi(\sqrt{N})}[\pi(N/p_k)-k+1]=((\sum_{k=1}^{\pi((N/2))}[\pi(N/p_k)]) + \pi (\sqrt{N}))/2$$
Simplifying LHS,
$$=\sum_{k=1}^{\pi(\sqrt{N})}[\pi(N/p_k)]-\sum_{k=1}^{\pi(\sqrt{N})}[k-1]$$
$$=\sum_{k=1}^{\pi(\sqrt{N})}[\pi(N/p_k)]-(\pi(\sqrt{N})(\pi(\sqrt{N})+1))/2+\pi(\sqrt{N})$$
$$=\sum_{k=1}^{\pi(\sqrt{N})}[\pi(N/p_k)]-((\pi(\sqrt{N}))^2-\pi(\sqrt{N})+2\pi(\sqrt{N}))/2$$
$$=(2(\sum_{k=1}^{\pi(\sqrt{N})}[\pi(N/p_k)])-(\pi(\sqrt{N}))^2+\pi(\sqrt{N}))/2$$
Equating LHS and RHS
$$(2(\sum_{k=1}^{\pi(\sqrt{N})}[\pi(N/p_k)])-(\pi(\sqrt{N}))^2+\pi(\sqrt{N}))/2=$$$$((\sum_{k=1}^{\pi((N/2))}[\pi(N/p_k)]) + \pi (\sqrt{N}))/2$$
Multiplying both sides by 2,
$$2(\sum_{k=1}^{\pi(\sqrt{N})}[\pi(N/p_k)])-(\pi(\sqrt{N}))^2+\pi(\sqrt{N})=$$$$(\sum_{k=1}^{\pi((N/2))}[\pi(N/p_k)]) + \pi (\sqrt{N})$$
Subtracting both sides by $(\pi(\sqrt{N})-(\pi(\sqrt{N}))^2)$
$$2(\sum_{k=1}^{\pi(\sqrt{N})}[\pi(N/p_k)])=(\sum_{k=1}^{\pi((N/2))}[\pi(N/p_k)])+(\pi(\sqrt{N}))^2$$
Subtracting both sides by $(\sum_{k=1}^{\pi((N/2))}[\pi(N/p_k)])$
$$2(\sum_{k=1}^{\pi(\sqrt{N})}[\pi(N/p_k)])-(\sum_{k=1}^{\pi((N/2))}[\pi(N/p_k)])=(\pi(\sqrt{N}))^2$$
which is equivalent to,
$$2(\sum_{k=1}^{\pi(\sqrt{N})}[\pi(N/p_k)])-\sum_{k=1}^{\pi(\sqrt{N})}[\pi(N/p_k)]-\sum_{k=\pi(\sqrt{N})+1}^{\pi(N/2)}[\pi(N/p_k)]$$
$$=(\pi(\sqrt{N}))^2$$
$$\implies\sum_{k=1}^{\pi(\sqrt{N})}[\pi(N/p_k)]-\sum_{k=\pi(\sqrt{N})+1}^{\pi(N/2)}[\pi(N/p_k)]=(\pi(\sqrt{N}))^2$$
\textbf{Example 2.}
\\For $N=25$, RHS=
$$(\pi(5))^2$$
$$=3^2=9$$
LHS=
$$\sum_{k=1}^{\pi(5)}[\pi(25/p_k)]-\sum_{k=\pi(5)+1}^{\pi(12.5)}[\pi(25/p_k)]$$
$$\pi (12.5)=5, \pi (8.33)=4, \pi (5)=3, \pi(3.57)=2, \pi(2.27)=1$$
$$=(5+4+3)-(2+1)$$
$$=12-3=9$$
Therefore, LHS=RHS\\
\section*{Conclusion}
In this research, a new formula for the semi-prime counting function was found. This has been shown to be a fairly straightforward method for determining the number of semi-prime numbers that do not exceed a given integer $N$. This yields an identity made up of many prime counting functions when compared to the previously known formula.\\
\section*{Forward looking}
Similar, to how two formulas for semi-prime counting functions were found, it might be possible to find more than one formula for k-almost prime numbers. A k-almost prime is a composite number that is the product of exactly $k$ (few of which might be equal) primes\cite{b2}. They then could be equated to find more such identities, which might give insights into the prime counting function itself.\\

\end{document}